\documentclass{article}

\usepackage[preprint]{neurips_2024}

\usepackage[utf8]{inputenc} 
\usepackage[T1]{fontenc}    
\usepackage{hyperref}       
\usepackage{url}            
\usepackage{booktabs}       
\usepackage{amsfonts}       
\usepackage{nicefrac}       
\usepackage{microtype}      
\usepackage{xcolor}         
\usepackage{amsmath, amssymb, amsthm}
\usepackage{enumitem}
\setlist[itemize]{leftmargin=*}

\title{Epistemic Confidence Statement via Extended Likelihood}

\author{%
  Youngjo~Lee 
  \\
  Department of Statistics \\
  Seoul National University \\
  \texttt{youngjo@snu.ac.kr} \\
}

\begin{document}

\maketitle

\begin{abstract}
Fisher's fiducial probability has recently attracted renewed attention under
the notion of epistemic confidence. Epistemic confidence statements can be
formulated through extended likelihoods, thereby clarifying several
long-standing controversies regarding its fiducial probability properties.
It establishes a direct connection between Fisher's epistemic notion of
confidence for observed data and Neyman's frequentist aleatory coverage
probability for future data, thereby enabling extension of epistemic
confidence statements for multidimensional parameters. We demonstrate how
higher-order asymptotic theory can be applied to refine the first-order
asymptotic epistemic confidence statements of the observed region, as a
direct consequence of extended likelihood property.
\end{abstract}


\section{Introduction}

There are two main interpretations of probability: aleatory probability such
as frequentist coverage probability, and epistemic probability such as the
Bayesian posterior. Fisher (1930) introduced fiducial probability as an
alternative to the Bayesian posterior that avoids the specification of a
prior, which enables epistemic confidence statements about parameters based
solely on observed data. However, fiducial inference has been the subject of
substantial controversy regarding its probabilistic foundations (Schweder
and Hjort, 2016). Recently, it has attracted renewed interest under the name
epistemic confidence (Pawtan et al., 2023;\ Lee, 2026) for hypothesis
testing and interval estimation. Schweder and Hjort (2016) applied
confidence distribution to interval estimation of scalar parameter. Pitman
(1957) demonstrated the inherent difficulty of extending confidence
(fiducial probability) to multi-dimensional parameters, since valid
confidence regions must then be restricted to particular forms. Several
alternatives have been proposed to address this limitation, including the
objective Bayesian approach using reference priors (Bernardo, 1979),
generalized fiducial inference (Hannig et al., 2016), and consonant belief
functions (Balch et al., 2019), each with its own advantages and domains of
applicability.

In satellite conjunction analysis, the assessed collision probability
approaches zero as observational noise increases, a phenomenon recognized as
a fundamental deficiency of probabilistic inference known as probability
dilution (Balch et al., 2019). Lee (2026) showed that epistemic confidence
can overcome probability dilution because it is not itself a probability,
whereas none of the above can fully resolve this deficiency. Recently,
confidence has been reinterpreted within an extended likelihood framework
(Pawitan \& Lee, 2021, 2024). Lee and Lee (2025) and Lee (2026) further
demonstrated several advantages of the epistemic confidence approach. It
resolves the induction problem of accepting a universal hypothesis based on
finite evidence while remaining falsifiable; it resolves Stein's (1959)
length problem while retaining exact coverage probability in a
high-dimensional setting; and it provides a coherent interpretation in
situations where frequentist coverage probability of an observed interval
becomes ambiguous.

Confidence does not satisfy Bayes' rule (Lee, 2026) and therefore cannot be
regarded as a Bayesian posterior. In this paper, it is reformulated as the
extended likelihood of a pivot, providing a coherent framework that links
epistemic confidence statements for an observed region with the frequentist
coverage probability of the corresponding region-generating procedure. This
formulation avoids long-standing controversies regarding its probabilistic
foundation. Moreover, it improves the accuracy of first-order asymptotic
epistemic confidence statements by incorporating recent advances in
higher-order frequentist asymptotic theory, while preserving Fisher's
original insight that epistemic confidence statements can be derived
directly from observed data. The proposed notion of confidence as an
extended likelihood thus retains both Neyman's (1937) frequentist coverage
interpretation and Fisher's (1930) epistemic interpretation.

\section{Fiducial probability and extended likelihood}

\subsection{Extended likelihood}

Let $\mathbf{Y}=(Y_{1},\cdots ,Y_{n})$ be a random vector from the joint
density function $f_{\mathbf{Y}^{\ast }}(\cdot |\boldsymbol{\theta}),$ where $%
\mathbf{Y}^{\ast }$ denotes an unobservable generic random vector having the
same joint density function as $\mathbf{Y}.$ Let $\mathbf{y}=(y_{1},\cdots
,y_{n})$ be the observation of $\mathbf{Y}$. Suppose that $\boldsymbol{\theta}%
=(\theta _{1},\cdots ,\theta _{q})$ is a fixed $q$-dimensional unknown
parameter. Given the observed data $\mathbf{Y}=\mathbf{y},$ $L(\boldsymbol{\theta};%
\mathbf{y})=f_{\mathbf{Y}^{\ast }}(\mathbf{y}|\boldsymbol{\theta})$ is the Fisher
(1921) likelihood for fixed unknown parameter$s$ $\boldsymbol{\theta}$, 
$\ell (\boldsymbol{\theta};\mathbf{y})=\log L(\boldsymbol{\theta};\mathbf{y})$ is the
log-likelihood and $\hat{\boldsymbol{\theta}}=\hat{\boldsymbol{\theta}}(\mathbf{y})$ is the
maximum likelihood estimator (MLE) of $\boldsymbol{\theta}$. In this paper, we use $%
a\equiv b$ to denote that $a$ is defined as $b$, semi-colon to separate
unknown quantities from observed data in the likelihood, and $F_{\mathbf{Y}%
^{\ast }}(\cdot |\boldsymbol{\theta})$ to denote the cumulative distribution
function for $\mathbf{Y}^{\ast }$. Suppose that $Y_{n+1}$ is a not-yet
observed future observation. Given observations $\mathbf{Y}=\mathbf{y},$ $%
L_{e}(Y_{n+1},\boldsymbol{\theta};\mathbf{y})\equiv f_{\mathbf{Y}^{\ast
},Y_{n+1}^{\ast }}(\mathbf{y},Y_{n+1}|\boldsymbol{\theta})$ is the extended
likelihood for prediction of an unobserved future observation $Y_{n+1}.$
When $Y_{n+1}$ is realized as $y_{n+1}$ but remains unobserved, as in the
case of latent variables or missing data. Then, the extended likelihood $%
L_{e}(y_{n+1},\boldsymbol{\theta};\mathbf{y})=f_{\mathbf{Y}^{\ast },Y_{n+1}^{\ast }}(%
\mathbf{y},y_{n+1}|\boldsymbol{\theta})$ can be used for inference of a realized
fixed unknown $y_{n+1}$ (Lee et al., 2017). For example, in animal breeding,
the breeding value $y_{n+1}$ of an already-born animal is realized as a
fixed unknown, whereas breeding value $Y_{n+1}$ of mating between two
potential parents is a random unknown. The joint density function $f_{%
\mathbf{Y}^{\ast },Y_{n+1}^{\ast }}(\mathbf{y},\cdot |\boldsymbol{\theta})$, which
provides the generating mechanism of the random unknown $Y_{n+1}$ and the
fixed unknown $y_{n+1}$, is assumed to be known up to a fixed unknown
parameter $\boldsymbol{\theta}$. No further distributional assumption of the fixed
unknown parameters $\boldsymbol{\theta}$ will be imposed in this paper. 

\subsection{Review of Fisher's confidence density}

We first review Fisher's (1930) confidence density, originally introduced as
the fiducial probability for a scalar parameter $\theta $. Given an
observation $\mathbf{Y}=\mathbf{y}$, let $s=s(\mathbf{y})$ be the realized
value of a scalar sufficient statistic, $S=s(\mathbf{Y})$ for $\theta $.
Suppose that $S$ is a continuous variable. Let the upper P-value of $S^{\ast
}=s(\mathbf{Y}^{\ast })$ be 
\begin{equation}
w=w(s,\theta )\equiv P(S^{\ast }\geq s|\theta )=1-F_{S^{\ast }}(s|\theta ),
\label{eq:CD}
\end{equation}
where $F_{S^{\ast }}(\cdot |\theta )$ denotes the cumulative distribution
function of $S^{\ast }$, with $\mathbf{Y}^{\ast }$ representing a generic
replicate drawn from the model parameterized by $\theta $. For a given $S=s,$
Fisher (1930) defined the confidence density and the likelihood as follows, 
\begin{equation}
c(\theta ;s)\equiv \frac{dP(S^{\ast }\geq s|\theta )}{d\theta }=\frac{dw}{%
d\theta }\quad \text{and}\quad L(\theta ;s)=\frac{dP(S^{\ast }\leq s|\theta )%
}{ds}=-\frac{dw}{ds}.  \label{eq:CL}
\end{equation}
For given data $\mathbf{Y}=\mathbf{y},$ they are the derivatives of a
realized pivot $w$ with respect to $\theta $ and $s$.

For subsequent notational development in the multi-dimensional parameter
case, consider the unobservable random pivot $W,$ which follows a uniform
distribution on $[0,1]$, 
$$
W\equiv w(S,\theta )=1-F_{S^{\ast }}(S|\theta ). 
$$
By construction, $W\sim $ Uniform$[0,1]$ under the sampling distribution
parameterized by $\theta $. The distinction between $S$ and its generic
version $S^{\ast }$ highlights that, in the frequentist approach, the
distribution function is defined with respect to repeated sampling through $%
S^{\ast }$, rather than the realized sample $S=s$. This distinction is
central to the development of epistemic confidence theory. Here, $%
w=1-F_{S^{\ast }}(s|\theta )$ is the realized value of the uniform pivot $W$%
, but it remains unobservable because $\theta $ is unknown. Thus, $w$
behaves like a fixed unknown parameter, although its generating mechanism is
known.

Suppose that $S^{\ast }=\hat{\theta}^{\ast }=\hat{\theta}(\mathbf{Y}^{\ast
}),$ and that $s=\hat{\theta}$ is sufficient. For a significance level $%
\alpha \ $with $0<\alpha <1$, let $g_{\alpha }(\cdot )$ be a strictly
increasing function; for example, $g_{\alpha }(x)=x+k$\ with inverse $%
g_{\alpha }^{-1}(x)=x-k.$ Consider the one-sided interval procedure $CI(%
\mathbf{Y}^{\ast })=[g_{\alpha }^{-1}(\hat{\theta}^{\ast }),\infty ).$
Define the unobservable binary random variable 
\begin{equation}
U^{\ast }=I(\theta \in CI(\mathbf{Y}^{\ast }))=I(\theta \geq g_{\alpha
}^{-1}(\hat{\theta}^{\ast }))=I(\hat{\theta}^{\ast }\leq g_{\alpha }(\theta
)).  \label{eq:UU}
\end{equation}
For any $\theta ,$ the extended likelihood of an unobserved event $U^{\ast
}=1$ is the coverage probability of the interval-generating procedure: 
$$
L_{{\mbox{e}}}(U^{\ast }=1)\equiv P(U^{\ast }=1|\theta )=P(\hat{\theta}%
^{\ast }\leq g_{\alpha }(\theta )|\theta )=1-\alpha . 
$$
This shows that the frequentist coverage probability is a property of the
generic random variable $\hat{\theta}^{\ast }$. Given the observed data $%
\mathbf{y}$, the generic interval-generating procedure $CI(\mathbf{Y}^{\ast
})$ is realized as the observed interval $CI(\mathbf{y})$ and the generic
random variable $U^{\ast }$ is realized as 
$$
u=I(\theta \in CI(\mathbf{y}))=I(\theta \geq g_{\alpha }^{-1}(\hat{\theta}%
))=I(\hat{\theta}\leq g_{\alpha }(\theta )), 
$$
which takes value 0 or 1 but remains unknown because $\theta $ is unknown.
Pawitan \& Lee (2021) showed that the extended likelihood $L_{\mbox{e}}(u;%
\hat{\theta})$ can be written as 
\begin{equation}
L_{{\mbox{e}}}(u=1;\hat{\theta})=C(\theta \in CI(\mathbf{y});\hat{\theta}%
)=\int_{\theta ^{\prime }\geq g_{\alpha }^{-1}(\hat{\theta})}c(\theta
^{\prime };\hat{\theta})d\theta ^{\prime }=1-\alpha ,
\end{equation}
which constitutes an epistemic confidence statement about the observed
interval. For the observed data, both $\theta $ and $\hat{\theta}$ are fixed
constants; hence, it cannot be a probability statement. The coverage
probability of the interval-generating procedure is referred to as the
confidence level. In this paper, we reserve the term confidence for the
epistemic interpretation of the observed interval $CI(\mathbf{y})$, while
using coverage probability to denote the frequentist confidence level of the
interval-generating procedure $CI(\mathbf{Y}^{\ast })$.

When $\hat{\theta}$\ is not sufficient, let $\mathbf{a}=a(\mathbf{y})\in 
\Bbb{R}^{n-q}$ be an ancillary statistic such that $(\hat{\theta},\mathbf{a}%
) $ is in one-to-one correspondence with $\mathbf{y}$ and is sufficient for $%
\theta \in \Bbb{R}^{q}$. Then the conditional density $f_{\widehat{\theta }%
^{\ast }}(\hat{\theta}|\mathbf{a},\theta )$ provides the likelihood based on
the full data information: 
$$
L(\theta ;\hat{\theta},\mathbf{a})\equiv f_{\hat{\theta}^{\ast }}(\hat{\theta%
}|\mathbf{a},\theta )\varpropto L(\theta ;\mathbf{y}). 
$$
Pawitan et al. (2023) showed that 
\begin{equation}
c(\theta ;\hat{\theta},\mathbf{a})=dP(\hat{\theta}^{\ast }\geq \hat{\theta}%
|\theta ,\mathbf{a})/d\theta  \label{eq:ancil}
\end{equation}
uses the full-data information with no relevant subset.

\section{Formulating confidence via extended likelihood}

Given data, suppose that we aim to make an epistemic confidence statement $%
C(\cdot )$\ of the observed region $CI(\mathbf{y})$\ for a $q$-dimensional $%
\boldsymbol{\theta}$. In this paper, $C(\cdot )$ is referred to as an epistemic
confidence statement if it satisfies the following confidence feature: 
\begin{equation}
C(\boldsymbol{\theta}\in CI(\mathbf{y});\mathbf{y})=P(\boldsymbol{\theta}\in CI(\mathbf{Y}%
^{\ast })|\boldsymbol{\theta}).  \label{eq:CF}
\end{equation}
Fisher (1958) insisted on the necessity of an epistemic interpretation 
$$
L_{\mbox{e}}(u=1;\mathbf{y})=C(\boldsymbol{\theta}\in CI(\mathbf{y});\mathbf{y})
$$
of the confidence statement for the observed region because scientists are
unlikely to repeat the same experiment many times. Neyman (1937) advocated
the frequentist aleatory interpretation 
$$
L_{\mbox{e}}(U^{\ast }=1)=P(\boldsymbol{\theta}\in CI(\mathbf{Y}^{\ast })|\boldsymbol{\theta}%
)\ 
$$
of the coverage probability for the region-generating procedure. Via the
confidence feature, it is possible to accommodate both the epistemic
confidence interpretation of the observed region and the frequentist
coverage probability statement for the region-generating procedure. Neyman
(1934) understood that the epistemic confidence statement (fiducial
probability) could serve as a tool for obtaining the coverage probability of
region-generating procedures. In contrast, via the confidence feature, the
frequentist coverage probability can also give an epistemic confidence
statement for the observed region. Thus, epistemic confidence and
frequentist aleatory probability play complementary roles in statistical
induction: see Lee (2026) for further discussion.

Suppose that we have a scalar pivot $V=v(\boldsymbol{\theta},\mathbf{Y}).$ Since
there is a one-to-one correspondence between $\mathbf{Y}$ and $(\hat{%
\boldsymbol{\theta}}(\mathbf{Y)},\mathbf{a})$, $V=v(\boldsymbol{\theta},\hat{\boldsymbol{\theta}}(%
\mathbf{Y)},\mathbf{a})$. Similarly, 
$$
W=w(\boldsymbol{\theta},\mathbf{Y})=w(\boldsymbol{\theta},\hat{\boldsymbol{\theta}}(\mathbf{Y)},%
\mathbf{a})=F_{V^{\ast }}(V), 
$$
is also a pivot and follows the uniform distribution. Its realized value is 
\begin{equation}
w=w(\boldsymbol{\theta},\mathbf{y})=P(V^{\ast }\leq v)=F_{V^{\ast }}(v)
\label{eq:pivot}
\end{equation}
where $v=v(\boldsymbol{\theta},\mathbf{y})$. Both $w$ and $v$ are fixed unknown
quantities, like parameters.

First, we show that epistemic confidence statements can be made solely based
on the distribution of generic random variables, $V^{\ast }=v(\boldsymbol{\theta},%
\mathbf{Y}^{\ast })=v(\boldsymbol{\theta},\hat{\boldsymbol{\theta}}(\mathbf{Y^{\ast })},%
\mathbf{a}^{\ast })$ and $W^{\ast }=w(\boldsymbol{\theta},\mathbf{Y}^{\ast })=w(%
\boldsymbol{\theta},\hat{\boldsymbol{\theta}}(\mathbf{Y^{\ast })},\mathbf{a}^{\ast })$.
Following \eqref{eq:CL}, we define the extended likelihoods as the
derivatives of a realized pivot\ with respect to the quantities of interest
as follows 
\begin{equation}
L_{{\mbox{e}}}(v;\mathbf{y})=\frac{dF_{V^{\ast }}(v)}{dv}=\frac{dw}{dv}%
=f_{V^{\ast }}(v)\quad \text{and}\quad L_{{\mbox{e}}}(V)=\frac{dF_{V^{\ast
}}(V)}{dV}=\frac{dw}{dV}=f_{V^{\ast }}(V).  \label{eq:EL}
\end{equation}
The extended likelihood $L_{\mbox{e}}(V)$ for a random $V=v(\boldsymbol{\theta},%
\mathbf{Y})$ behaves like an aleatory probability of a generic random
variable $V^{\ast }$, whereas the extended likelihood $L_{\mbox{e}}(v;%
\mathbf{y})$ for a fixed unknown $v=v(\boldsymbol{\theta},\mathbf{y})$ behaves like
the likelihood of a fixed unknown parameter $v$.

In this paper, we consider an observed region of the form $CI(\mathbf{y})=\{%
\boldsymbol{\theta}:v\leq b\}$ and a region-generating procedure $CI(\mathbf{Y}%
^{\ast })=\{\boldsymbol{\theta}:V^{\ast }\leq b)$. Since the statements $\boldsymbol{\theta}%
\in CI(\mathbf{y})$ and $v\leq b$\ are equivalent and both $\boldsymbol{\theta}$ and 
$v$ are fixed unknown quantities, we define the confidence statement 
$$
C(\boldsymbol{\theta}\in CI(\mathbf{y});\mathbf{y})\equiv C(v\leq b;\mathbf{y}).
$$
This is consistent with Fisher's (1959) fiducial argument based on the $t$%
-pivot discussed in the next section. The controversies stem from treating $%
C(\boldsymbol{\theta}\in CI(\mathbf{y});\mathbf{y})$ as a probability statement,
although after the data have been observed, both the unknown parameter $%
\boldsymbol{\theta}$ and the realized confidence set $CI(\mathbf{y})$ are fixed
quantities. 

Let 
$$
L_{\mbox{e}}(v;\mathbf{y})=f_{V^{\ast }}(v)\equiv c(v;\mathbf{y}).
$$
From (\ref{eq:EL}), the confidence feature (\ref{eq:CF}) follows immediately
because the epistemic confidence statement for the observed region $CI(%
\mathbf{y}),$ 
$$
C(\boldsymbol{\theta}\in CI(\mathbf{y});\mathbf{y})=C(v\leq b;\mathbf{y})=F_{V^{\ast
}}(b),
$$
and the aleatory coverage probability of the region-generating procedure $CI(%
\mathbf{Y}^{\ast }),$%
$$
P(\boldsymbol{\theta}\in CI(\mathbf{Y}^{\ast })|\boldsymbol{\theta})=P(V^{\ast }\leq
b)=F_{V^{\ast }}(b),
$$
are the same. Thus, the confidence feature is simply a consequence of
extended likelihood property 
$$
L_{\mbox{e}}(v;\mathbf{y})=c(v;\mathbf{y})=f_{V^{\ast }}(v)=L_{\mbox{e}%
}(V^{\ast }=v).
$$
For the epistemic confidence statement $C(\boldsymbol{\theta}\in CI(\mathbf{y});%
\mathbf{y})$, there is no need to define a confidence density $c(\boldsymbol{\theta};%
\mathbf{y})$\ for $\boldsymbol{\theta}$; it suffices to specify the density of the
generic random variable $V^{\ast }$. Thus, there is no need to restrict
valid confidence regions to particular forms (Pitman, 1957) in
mult-dimensional cases. When $V$ is the likelihood ratio statistic, $P(%
\boldsymbol{\theta}\in CI(\mathbf{Y}^{\ast })|\boldsymbol{\theta})=P(V^{\ast }\leq b|%
\boldsymbol{\theta})$ is a standard way of computing the coverage probability.
Expressing confidence statements through extended likelihood enables us to
reconcile Fisher's epistemic interpretation with Neyman's aleatory
interpretation, even for multidimensional parameters.

When $q=1$, from \eqref{eq:CL}, Fisher's confidence density of $\theta $ and
the density of $\hat{\theta}$ can be derived from the extended likelihood of 
$v$ through the transformation rule: 
\begin{eqnarray}
c(\theta ;\mathbf{y}) &=&L_{{\mbox{e}}}(v;\mathbf{y})\left| \frac{dv}{%
d\theta }\right| =f_{V^{\ast }}(v)\left| \frac{dv}{d\theta }\right| ,
\label{eq:tt} \\
f_{\hat{\theta}^{\ast }}(\hat{\theta}|\theta ) &=&f_{V^{\ast }}(v)\left| 
\frac{dv}{d\hat{\theta}}\right| =L_{{\mbox{e}}}(v;\mathbf{y})\left| \frac{dv%
}{d\hat{\theta}}\right| ,  \nonumber
\end{eqnarray}
where $f_{\hat{\theta}^{\ast }}(\hat{\theta}|\theta )$ is the density of $%
\hat{\theta}^{\ast }.$\ Thus, this extended likelihood formulation
accommodates Fisher's original confidence density for a scalar parameter $%
\theta .$

Example 1: Let $q=1.$ Suppose that $\hat{\theta}\sim N(\theta ,1).$ Then, 
$$
V^{\ast }=\hat{\theta}^{\ast }-\theta \sim N(0,1)
$$
is a pivot with its realized value 
$$
v=\hat{\theta}-\theta .
$$
The extended likelihood $L_{e}(v;\mathbf{y})=f_{V^{\ast }}(v)=\exp
(-v^{2}/2)/\sqrt{2\pi }=\phi (v)=c(v;\mathbf{y})$ where $\phi (\cdot )$
denotes the standard normal density function\ and $c(\theta ;\mathbf{y}%
)=f_{V^{\ast }}(v)|dv/d\theta |=f_{V^{\ast }}(v)$ with $|dv/d\theta |=1.$
Since $P(\hat{\theta}^{\ast }\geq \hat{\theta}|\theta )=1-\Phi (\hat{\theta}%
-\theta )$ with $\Phi (\cdot )$ being the cumulative distribution function
of standard normal distribution$,$ the extended likelihood provides Fisher's
confidence density (fiducial probability) of $\theta ,$ 
$$
c(\theta ;\mathbf{y})=dP(\hat{\theta}^{\ast }\geq \hat{\theta}|\theta
)/d\theta =d\{1-\Phi (\hat{\theta}-\theta ))/d\theta =\phi (\hat{\theta}%
-\theta )=\phi (v)=f_{V^{\ast }}(v),
$$
and the density function of the MLE $\hat{\theta}$, 
$$
f_{\hat{\theta}^{\ast }}(\hat{\theta}|\theta )=f_{V^{\ast }}(v)|dv/d\hat{%
\theta}|=f_{V^{\ast }}(v),
$$
so both $c(\theta ;\mathbf{y})$\ and $f_{\hat{\theta}^{\ast }}(\hat{\theta}%
|\theta )$ can be derived directly from the extended likelihood $L_{e}(v;%
\mathbf{y})=c(v;\mathbf{y}).$ Thus, the confidence feature (\ref{eq:CF})
follows immediately because 
$$
C(\boldsymbol{\theta}\in CI(\mathbf{y});\mathbf{y})=C(v\leq b;\mathbf{y}%
)=\int_{v\leq b}f(v)dv=\Phi (b)
$$
and 
$$
P(\boldsymbol{\theta}\in CI(\mathbf{Y}^{\ast })|\boldsymbol{\theta})=P(V^{\ast }\leq b)=\Phi
(b).
$$
Fisher's fiducial probability serves as an instrument enabling an epistemic
interpretation of the observed region. Extended likelihood gives a natural
way to establish confidence feature by connecting Fisher's epistemic
statement based on $v$\ with Neyman's frequentist coverage probability based
on $V^{\ast }$ without giving rise to the usual controversies. As we shall
show, it also provides an immediate extension of Fisher's epistemic
confidence statements to multi-dimensional parameters.


\subsection{Epistemic confidence statement for multi-dimensional parameters}

With $\mathbf{j}(\boldsymbol{\theta})=-\partial ^{2}\ell (\boldsymbol{\theta};\mathbf{y}%
)/\partial \boldsymbol{\theta}\partial \boldsymbol{\theta}^{T}$, let $\hat{\mathbf{j}}=%
\mathbf{j}(\hat{\boldsymbol{\theta}})$ be the observed Fisher information. It is
immediate that the deviance, twice the log-likelihood ratio, 
\begin{eqnarray}
d=d(\boldsymbol{\theta},\mathbf{y}) &\equiv &2\{\ell (\hat{\boldsymbol{\theta}};\mathbf{y}%
)-\ell (\boldsymbol{\theta};\mathbf{y})\}  \label{eq:logLR} \\
&=&(\boldsymbol{\theta}-\hat{\boldsymbol{\theta}})^{T}\hat{\mathbf{j}}(\boldsymbol{\theta}-\hat{%
\boldsymbol{\theta}})\left\{ 1+O_{p}(n^{-1/2})\right\} ,  \label{eq:RL}
\end{eqnarray}
is an asymptotically $\chi _{q}^{2}$\ pivot. Let $L_{{\ \mbox{n}}}(%
\boldsymbol{\theta};\mathbf{y})=L(\boldsymbol{\theta};\mathbf{y})/\int L(\boldsymbol{\theta}^{\prime
};\mathbf{y})d\boldsymbol{\theta}^{\prime }$\ be the normalized likelihood. The
Laplace approximation to the integrated likelihood 
$$
\int L(\boldsymbol{\theta}^{\prime };\mathbf{y})d\boldsymbol{\theta}^{\prime }=|2\pi \mathbf{%
\hat{\mathbf{j}}}^{-1}|^{1/2}L(\hat{\boldsymbol{\theta}};\mathbf{y})\left\{
1+O_{p}(n^{-1/2})\right\} 
$$
yields the normalized likelihood 
\begin{equation}
L_{{\ \mbox{n}}}(\boldsymbol{\theta};\mathbf{y})=|2\pi \mathbf{\hat{\mathbf{j}}}%
^{-1}|^{-1/2}\exp \big\{-(\boldsymbol{\theta}-\hat{\boldsymbol{\theta}})^{T}\mathbf{\hat{%
\mathbf{j}}}(\boldsymbol{\theta}-\hat{\boldsymbol{\theta}})/2\big\}\left\{
1+O_{p}(n^{-1/2})\right\} ,  \label{eq:NL}
\end{equation}
leading to another asymptotic normal pivot, $\mathbf{\hat{\mathbf{j}}}^{1/2}(%
\boldsymbol{\theta}-\hat{\boldsymbol{\theta}})$. Thus, by taking either $v=d=2\{\ell (\hat{%
\boldsymbol{\theta}};\mathbf{y})-\ell (\boldsymbol{\theta};\mathbf{y})\}$ or $v=\mathbf{\hat{%
\mathbf{j}}}^{1/2}(\boldsymbol{\theta}-\hat{\boldsymbol{\theta}})$, the confidence feature
therefore holds asymptotically under either the first-order chi-squared
approximation to the likelihood ratio in (\ref{eq:RL}) or the first-order
normal approximation to the MLE in (\ref{eq:NL}). The focus of this paper is
on how to obtain exact or higher-order asymptotic epistemic confidence
statements for an observed region. Since the confidence feature generally
holds only asymptotically for discrete data, we therefore concentrate on
continuous data where exact confidence could be attained.

Lee (2026) showed that under a particular scale $\boldsymbol{\eta}=\boldsymbol{\eta}(%
\boldsymbol{\theta},\mathbf{y}),$ the confidence density can be written as the
modified normalized likelihood 
\begin{equation}
c(\boldsymbol{\theta};\mathbf{y})=\frac{c_{0}(\boldsymbol{\theta};\mathbf{y})L(\boldsymbol{\theta};%
\mathbf{y})}{\int c_{0}(\boldsymbol{\theta}^{\prime };\mathbf{y})L(\boldsymbol{\theta}%
^{\prime };\mathbf{y})d\boldsymbol{\theta}^{\prime }}=|d\boldsymbol{\eta}/d\boldsymbol{\theta}|L_{{%
\mbox{n}}}(\boldsymbol{\eta};\mathbf{y}),  \label{eq:IM}
\end{equation}
where $L_{{\mbox{n}}}(\boldsymbol{\eta};\mathbf{y})$ is the normalized likelihood on
the $\boldsymbol{\eta}$\ scale and $c_{0}(\boldsymbol{\theta};\mathbf{y})=|d\boldsymbol{\eta}/d%
\boldsymbol{\theta}|$ is the Jacobian term, which may depend on the data. In Example
1, where $\hat{\theta}\sim N(\theta ,1)$, the pivot is $v=\hat{\theta}%
-\theta $ and $L_{\mbox{e}}(v;\hat{\theta})=\exp (-v^{2}/2)/\sqrt{2\pi }$.
Thus, the confidence density is precisely the modified normalized likelihood
at the scale $\eta =\theta .$ However, when $\hat{\theta}-\theta $ does not
follow normal distribution, an alternative scale $\eta $\ could improve the
first-order confidence statement based on normal approximation of $\hat{%
\theta}$\ in finite sample.

Let $\boldsymbol{\theta}=(\boldsymbol{\alpha},\boldsymbol{\gamma})$ with $\boldsymbol{\alpha}\in \Bbb{R}^{s}$
being parameters of interest and $\boldsymbol{\gamma}\in \Bbb{R}^{q-s}$ nuisance
parameters. The profile likelihood of $\boldsymbol{\alpha}$ is 
$$
L_{\mp }(\boldsymbol{\alpha};\mathbf{y})=L(\boldsymbol{\alpha},\tilde{\boldsymbol{\gamma}}(%
\boldsymbol{\alpha});\mathbf{y}), 
$$
where $\tilde{\boldsymbol{\gamma}}(\boldsymbol{\alpha})$ is the profile MLE of $\boldsymbol{\gamma}$
for a given $\boldsymbol{\alpha}.$ The profile deviance 
$$
d_{\mp }=2\{\ell _{\mp }(\hat{\boldsymbol{\alpha}};\mathbf{y})-\ell _{\mp }(%
\boldsymbol{\alpha};\mathbf{y})\} 
$$
gives an asymptotic $\chi _{s}^{2}$\ pivot. Thus, we can make a first-order
asymptotic epistemic confidence statement for $\boldsymbol{\alpha}$ using $d_{\mp}$
(Schweder \& Hjort, 2016).

We first investigate how to construct confidence statements when pivots are
explicitly available. Consider a normal regression model, 
$$
\mathbf{Y}=\mathbf{X\boldsymbol{\beta}}+\mathbf{e}, 
$$
where $\mathbf{X}$ is an $n\times p$ matrix and $\mathbf{e}\sim N(\mathbf{0}%
,\phi \mathbf{I}_{n}).$ Here $q=p+1.$ Let $\mathbf{\hat{\boldsymbol{\beta}}}=(%
\mathbf{X}^{T}\mathbf{X})^{-1}\mathbf{X}^{T}\mathbf{y}$ and 
$$
\hat{\phi}_{M}=(\mathbf{y}-\mathbf{X\hat{\boldsymbol{\beta}}})^{T}(\mathbf{y}-%
\mathbf{X\hat{\boldsymbol{\beta}}})/(n-p). 
$$

Suppose that we want to construct a confidence statement for $\boldsymbol{\alpha}%
=\phi $ with $s=1$. Let 
\begin{equation}
v=v(\phi {,\mathbf{y})}=(n-p)\hat{\phi}_{M}/\phi .  \label{eq:LRs}
\end{equation}
Then, $V^{\ast }=v(\phi \mathbf{,Y}^{\ast }\mathbf{)}$\ is a $\chi
_{n-p}^{2} $\ pivot to give 
$$
P(V^{\ast }\leq v)=F_{V^{\ast }}(v)=\gamma ((n-p)/2),v/2)/\Gamma ((n-p)/2), 
$$
where $\gamma (k,x)=\int_{0}^{x}t^{k-1}e^{-t}dt$ is the incomplete gamma
function$.$ Thus, 
$$
L_{{\mbox{e}}}(v;\mathbf{y})=c(v;\mathbf{y})=(v/2)^{(n-p)/2}v^{-1}\exp
(-v/2)/\Gamma ((n-p)/2)=f_{V^{\ast }}(v). 
$$
By \eqref{eq:LRs}, $dv/d\phi =-v/\phi $, 
$$
c(\phi ;\mathbf{y})=c(v;\mathbf{y})|dv/d\phi |=(v/\phi )f_{V^{\ast
}}(v)=\phi ^{-1}(v/2)^{(n-p)/2}\exp (-v/2)/\Gamma ((n-p)/2). 
$$
By \eqref{eq:LRs}, $dv/d\hat{\phi}_{M}=(n-p)/\phi ,$ 
$$
f_{\hat{\phi}_{M}^{\ast }}(\hat{\phi}_{M}|\phi )=\{(n-p)/\phi \}f_{V^{\ast
}}(v)=\hat{\phi}_{M}(v/2)^{(n-p)/2}\exp (-v/2)/\Gamma ((n-p)/2), 
$$
where $f_{\hat{\phi}_{M}^{\ast }}(\hat{\phi}_{M}|\phi )$ is the density of $%
\hat{\phi}_{M}^{\ast }$ with parameter $\phi $. These expressions imply that 
$$
\phi c(\phi ;\mathbf{y})=\hat{\phi}_{M}^{-1}f_{\hat{\phi}_{M}^{\ast }}(\hat{%
\phi}_{M}|\phi ). 
$$
Let $CI(\mathbf{y)}=\{\phi :v\leq b\}$ and $CI(\mathbf{Y}^{\ast }\mathbf{)}%
=\{\phi :V^{\ast }\leq b)$. Then, the confidence feature is immediate
because 
$$
C(\phi \in CI(\mathbf{y)})\equiv C(v\leq b;\mathbf{y})=F_{V^{\ast }}(b)\text{
\ and \ }P(\phi \in CI(\mathbf{Y}^{\ast }\mathbf{)}|\phi )=F_{V^{\ast }}(b), 
$$
where $C(v\leq b;\mathbf{y})$ is an epistemic confidence statement about
fixed unknown $v$ and $P(V^{\ast }\leq b)$ is the coverage probability of
the region-generating procedure$.$ By \eqref{eq:LRs}, $dv/d\phi =-v/\phi $
and $v\varpropto \hat{\phi}_{M}/\phi ,$ $L(\phi ;\mathbf{y})=f_{\hat{\phi}%
_{M}^{\ast }}(\hat{\phi}_{M}|\phi )=f_{V^{\ast }}(v)|dv/d\hat{\phi}%
_{M}|\varpropto f_{V^{\ast }}(v)/\phi $\ gives 
$$
c(\phi ;\mathbf{y})=c(v;\mathbf{y})|dv/d\phi |=(v/\phi )f_{V^{\ast
}}(v)\varpropto (\hat{\phi}_{M}/\phi )L(\phi ;\mathbf{y}). 
$$
Thus, with $\eta =\hat{\phi}_{M}\log \phi $, the Fisher confidence density
for $\phi $, $c(\phi ;\mathbf{y})$, is the modified normalized likelihood at
the scale $\eta =\hat{\phi}_{M}\log \phi $ in (\ref{eq:IM}).

Let $\lambda =\mathbf{u}^{T}\boldsymbol{\beta}$ be the parameter of interest with $%
s=1,$ $\hat{\lambda}=\mathbf{u}^{T}\hat{\boldsymbol{\beta}}$ and $k=\mathbf{u}^{T}(%
\mathbf{X}^{T}\mathbf{X})^{-1}\mathbf{u}.$ Then, 
$$
v=t=(\hat{\lambda}-\lambda )/(k\hat{\phi}_{M})^{1/2}\sim t_{n-p},
$$
where $t_{m}$ denotes a $t$-distribution with $m$ degrees of freedom. Since $%
|dv/d\lambda |=1/(k\hat{\phi}_{M})^{1/2}$, we have $c(\lambda
;y)=c(v;y)|dv/d\lambda |.$ In this example, let $CI(\mathbf{y)}=\{\lambda
:v\leq b\}$ and $CI(\mathbf{Y}^{\ast }\mathbf{)}=\{\lambda :V^{\ast }\leq b).
$ Then, the confidence feature is immediate because 
$$
C(\lambda \in CI(\mathbf{y)};\mathbf{y})\equiv C(v\leq b;\mathbf{y}%
)=F_{V^{\ast }}(b)\text{ \ and \ }P(\lambda \in CI(\mathbf{Y}^{\ast
})|\lambda )=P(V^{\ast }\leq b)=F_{V^{\ast }}(b).
$$
Fisher (1957) denoted the confidence statement 
$$
C((\hat{\lambda}-\lambda )/(k\hat{\phi}_{M})^{1/2}\leq b;y)=C(t\leq b;y)
$$
by the probability statement 
$$
P(t\leq b),
$$
which led to unnecessary controversies because both $t\ $and $b$ are
constants once the data are observed. With $V^{\ast }=(\hat{\lambda}^{\ast
}-\lambda )/(k\hat{\phi}_{M}^{\ast })^{1/2}=t^{\ast },$ the corresponding
probability statement 
$$
P(V^{\ast }\leq b)=P\{(\hat{\lambda}^{\ast }-\lambda )/(k\hat{\phi}%
_{M}^{\ast })^{1/2}\leq b\}=P(t^{\ast }\leq b),
$$
is Neyman's coverage probability. To circumvent these controversies, we
denote the epistemic confidence statement for the observed region by 
$$
C(\lambda \in CI(\mathbf{y)};\mathbf{y})\equiv C(v\leq b;\mathbf{y})=C(t\leq
b;\mathbf{y})
$$
and the coverage probability statement for the region-generating procedure
by 
$$
P(\lambda \in CI(\mathbf{Y}^{\ast })|\lambda )=P(V^{\ast }\leq b)=P(t^{\ast
}\leq b).
$$
Extended likelihoods $L_{e}\left( V\right) $ and $L_{e}\left( v;\mathbf{y}%
\right) $ provide the frequentist coverage probability of the
region-generating procedure and the epistemic confidence statement for the
observed region, respectively. Fisher's (1957) fiducial probability
statement is more properly interpreted as epistemic confidence statement $%
C(t\leq b;y)$ based on the pivot $t$.

Now an epistemic confidence statement can be made straightforwardly for
multi-dimensional parameters $\boldsymbol{\beta}\in \Bbb{R}^{p}$ with $s=p.$ Let $v=(%
\hat{\boldsymbol{\beta}}-\boldsymbol{\beta})^{T}\mathbf{X}^{T}\mathbf{X}(\hat{\boldsymbol{\beta}}-%
\boldsymbol{\beta})/(p\hat{\phi}_{M})$ be a pivot, which follows an $F$%
-distribution. Let $CI(\mathbf{y)}=\{\eta :v\leq b\}$ and $CI(\mathbf{Y}%
^{\ast })=\{\boldsymbol{\beta}:V^{\ast }\leq b\}$ to show the confidence feature (%
\ref{eq:CF}) for $\boldsymbol{\beta}$ 
$$
C(\boldsymbol{\beta}\in CI(\mathbf{y);y})\equiv C(v\leq b\mathbf{;y})=\int_{v\leq
b}f_{_{V^{\ast }}}(v)dv=F_{_{V^{\ast }}}(b)\text{ \ } 
$$
and 
$$
P(\boldsymbol{\beta}\in CI(\mathbf{Y}^{\ast })|\boldsymbol{\beta})=P\{V^{\ast }\leq
b\}=F_{_{V^{\ast }}}(b). 
$$
Due to this confidence feature, epistemic confidence can give coverage
probability of a region-generating procedure, as Neyman (1934) understood.
Furthermore, the coverage probabilities can also provide epistemic
confidence statements.

\subsection{Improving epistemic confidence statements via higher-order
approximations}

We shall show how to improve the first-order asymptotic epistemic confidence
statement of observed region. An alternative approximation to the normalized
likelihood in (\ref{eq:NL}) is to use (\ref{eq:logLR}) directly instead of (%
\ref{eq:RL}) 
$$
L_{{\mbox{n}}}(\boldsymbol{\theta};\mathbf{y})=|\mathbf{\hat{\mathbf{j}}}/(2\pi
)|^{1/2}\exp \big\{-[\ell (\hat{\boldsymbol{\theta}};\mathbf{y})-\ell (\boldsymbol{\theta};%
\mathbf{y})]\big\}(1+O_{p}(n^{-1/2})),
$$
whose renormalization leads to Barndorff-Nielsen's (1983) higher-order
approximation to the density function of the MLE $\hat{\boldsymbol{\theta}}$: 
\begin{equation}
f_{\hat{\boldsymbol{\theta}}^{\ast }}(\hat{\boldsymbol{\theta}}|\mathbf{a},\boldsymbol{\theta})=b|%
\mathbf{\hat{\mathbf{j}}}/(2\pi )|^{1/2}\exp \big\{-[\ell (\hat{\boldsymbol{\theta}};%
\mathbf{y})-\ell (\boldsymbol{\theta};\mathbf{y})]\big\}(1+O_{p}(n^{-1})),
\label{eq:high}
\end{equation}
where $b$ is the normalizing constant ensuring $\int f_{\hat{\boldsymbol{\theta}}%
^{\ast }}(\hat{\boldsymbol{\theta}}|\mathbf{a},\boldsymbol{\theta})d\hat{\boldsymbol{\theta}}=1$.
The higher-order asymptotic theory is an essential component of frequentist
inference (Reid, 2003)$.$ There is a long history of extending higher-order
approximations in (\ref{eq:high}) to improve the accuracy of coverage
probability computation in finite samples, rather than relying on the
first-order normal approximation in (\ref{eq:NL}) or the chi-squared
limiting distribution in (\ref{eq:RL}).

Barndorff-Nielsen's (1983) higher-order approximation does not require
explicit specification of the ancillary statistic $\mathbf{a}$. Similarly,
in this section we assume only the existence of a pivot with either a
chi-squared or a normal distribution, without requiring its explicit form.
If a pivot is explicitly available, then the extended likelihood in (\ref
{eq:EL}) provides an exact epistemic confidence statement and also enables
exact computation of the corresponding coverage probability, as shown in the
previous section. However, if the pivot is not explicitly available, as we
shall demonstrate, higher-order asymptotic can improve first-order
asymptotic epistemic confidence statements. Higher-order asymptotic theories
have been extended to the parameter of interest via the profile deviance.

For illustration, we consider a gamma regression model with the
log-likelihood 
$$
\ell (\boldsymbol{\beta},\varphi ;\mathbf{y})=-\varphi \sum \{(y_{i}-\mu _{i})/\mu
_{i}-\log (y_{i}/\mu _{i})\}-n\{\log \Gamma (\varphi )-\varphi \log (\varphi
)+\varphi \}, 
$$
where $\varphi =1/\phi $ is the precision parameter, $\mu _{i}=\exp (\mathbf{%
x}_{i}^{T}\boldsymbol{\beta})$ and $\Gamma (\cdot )$ is the gamma function$.$

Suppose that we want to make a confidence statement of $\mathbf{\boldsymbol{\alpha}=}%
\varphi $ with $s=1.$\ Let $d_{\mp }=2\{\ell _{\mp }(\hat{\varphi};\mathbf{y}%
)-\ell _{\mp }(\varphi ;\mathbf{y})\}$ be the profile deviance. Improvement
of the first-order normal approximation\ to $z_{\mp }=\mathrm{sign}(\hat{%
\varphi}-\varphi )d_{\mp }^{1/2}$ can be made by using the higher-order
approximation 
\begin{equation}
P(Z^{\ast }\leq z|\varphi )=F_{Z^{\ast }}(z|\varphi )=\Phi
(z)\{1+O(n^{-1})\},  \label{eq:ZZ}
\end{equation}
where $Z^{\ast }$ is the unknown generic standard normal pivot, 
\begin{equation}
z=z(\varphi ,\mathbf{y})=z_{\mp }+(1/z_{\mp })\log (m/z_{\mp }),
\label{eq:ZP}
\end{equation}
and expressions for $m$ have been proposed, for example by Fraser (1991),
Skovgaard (2001), Tang and Reid (2025). Taking $m=z_{\mp }$ results in $%
z=z_{\mp }$, thereby yielding the first-order asymptotic coverage
probability of the log-likelihood ratio interval. It is important to see the
analogy between the extended likelihood (\ref{eq:pivot}) with $V^{\ast
}=Z^{\ast }$ and the higher-order asymptotic (\ref{eq:ZP}) with $%
v=z=z(\varphi ,\mathbf{y}).$ Higher-order asymptotics have been developed to
compute the frequentist coverage probability at each value of $\varphi $. If
the expression of $Z^{\ast }$\ is known explicitly, the approximation is
exact $P(Z^{\ast }\leq z)=\Phi (z)$ without depending on $\varphi $. Then,
confidence feature holds immediately as in the previous section; the
coverage probability should remain the same\ as the stated coverage
probability for all $\varphi $. However, when the approximation is not
exact, for each value of $\varphi ,$ (\ref{eq:ZZ}) yields approximate
coverage probabilities that are substantially closer to the nominal level
than those based on the first-order normal approximation; see the simulation
studies of Skovgaard (1991) and Tang and Reid (2025). Pawitan et al. (2023)
studied the computation of epistemic confidence statements when $s=1$\ and
showed that their computed confidence statement $C(\varphi \in CI(\mathbf{%
y);y})=C(\varphi \in CI(\mathbf{y);\hat{\boldsymbol{\theta}},\mathbf{a)}}$ for given
data is essentially identical to the coverage probability $P(\varphi \in CI(%
\mathbf{Y}^{\ast })|\varphi =\hat{\varphi})$ evaluated at the MLE $\varphi =%
\hat{\varphi}$ using the higher-order asymptotic approximation of Fraser
(1991), which, in this paper, we call the maximum likelihood epistemic
confidence level.

Consider the confidence statement of $\mathbf{\boldsymbol{\alpha}=}\boldsymbol{\beta}\in 
\Bbb{R}^{p}$ with $s=p$. To improve the first-order asymptotic $\chi _{p}^{2}
$ approximaton to $d_{\mp }=2\{\ell _{\mp }(\hat{\boldsymbol{\beta}})-\ell _{\mp }(%
\boldsymbol{\beta})\},$ Skovgaard (2001) proposed a higher-order approximation 
$$
P(D^{\ast }\leq d|\boldsymbol{\beta})=F_{D^{\ast }}(d)=\frac{\gamma (p/2,d/2)}{%
\Gamma (p/2)}\{1+O(n^{-1})\},
$$
where $D^{\ast }$ is the generic $\chi _{p}^{2}$ random variable, $d=d(%
\boldsymbol{\beta},\mathbf{y})=d_{\mp }-2\log (m),$ an expression of $m=m(\boldsymbol{\beta},%
\mathbf{y})$ is given in Skovgaard (2001) and $\gamma (k,x)$ is the
incomplete gamma function. In this paper, $d$ is a higher-order
approximation to the realized value of the unknown pivot $D^{\ast }$\ that
yields an accurate computation of the coverage probability at each value of $%
\boldsymbol{\beta}$. Let $CI(\mathbf{y)}=\{\boldsymbol{\beta}:d\leq k\}$ and $CI(\mathbf{Y}%
^{\ast })=\{\boldsymbol{\beta}:D^{\ast }\leq k\}.$ From (\ref{eq:CF}), if the
approximation $P(D^{\ast }\leq d)=\gamma (p/2,d/2)/\Gamma (p/2)$\ is exact,
the previous section shows that the confidence feature for the
multi-dimensional parameter $\boldsymbol{\beta}$ is immediate because they are the
same quantity: 
$$
C(\boldsymbol{\beta}\in CI(\mathbf{y)})\equiv C(d\leq k;\mathbf{y})=F_{D^{\ast }}(k)%
\text{ \ and \ }P(\boldsymbol{\beta}\in CI(\mathbf{Y}^{\ast }))=F_{D^{\ast }}(k).
$$
If the explicit expression of the pivot is not known, as the maximum
likelihood epistemic confidence statement $C(\boldsymbol{\beta}\in CI(\mathbf{y);%
\hat{\boldsymbol{\theta}},\mathbf{a}})$ for the observed region $CI(\mathbf{y)}$ we
use the frequentist coverage probability $P(\boldsymbol{\beta}\in CI(\mathbf{Y}%
^{\ast })|\boldsymbol{\beta})$ based on a higher-order asymptotic approximation,
evaluated at $\boldsymbol{\beta}=\hat{\boldsymbol{\beta}}$. Therefore, the higher-order
approximation to Neyman's coverage probability refines Fisher's first-order
asymptotic epistemic confidence statement based on chi-square approximation
to the profile deviance (Schweder \& Hjort, 2016).

\section{Conclusion}

Lee and Lee (2025) and Lee (2026) articulated both conceptual and practical
advantages of the confidence approach. The extended likelihood provides a
coherent framework for resolving long-standing controversies surrounding
fiducial probability while preserving Fisher's original insight that
epistemic confidence statements can be extracted from observed data and
maintaining Neyman's frequentist coverage interpretation. Consequently, a
frequentist region-generating procedure for multidimensional parameters with
a specified coverage probability admits a corresponding epistemic confidence
interpretation of the observed region. Higher-order approximations can
improve the accuracy of first-order asymptotic epistemic confidence
statements. The epistemic and aleatory interpretations of probability can be
reconciled through framework of extended likelihood, emphasizing the synergy
rather than the conflict between epistemic and aleatory interpretations.

\section*{References}

\begin{itemize}
\item  {Balch, M. S., Martin, R., \& Ferson, S.} (2019). Satellite
conjunction analysis and the false confidence theorem. \textit{Proceedings
of the Royal Society A}, \textbf{475}(2227), 20180565.

\item  {Barndorff-Nielsen, O.} (1983). On a formula for the distribution of
the MLE. \textit{Biometrika}, \textbf{70}(2), 343--365.

\item  {Bernardo, J. M.} (1979). Reference posterior distributions for
Bayesian inference. \textit{Journal of the Royal Statistical Society Series B%
}, \textbf{41}(2), 113--128.

\item  {Fisher, R. A. } (1921). On the probable error' of a coefficient of
correlation deduced from a small sample. \textit{Metron}, \textbf{1}, 1-32.

\item  {Fisher, R. A.} (1930). Inverse probability. \textit{Mathematical
Proceedings of the Cambridge Philosophical Society}, \textbf{26}(4),
528--535.

\item  {Fisher, R. A. } (1958). The Nature of Probability. 
\textit{Centennial Review of Arts and Science}, \textbf{2}, 261-274.

\item  {Fisher, R. A. } (1959). Mathematical probability in the natural science.
\textit{Techmetrics}, \textbf{1}, 21-29.

\item  {Fraser, D. A. S.} (1991). Statistical inference: Likelihood to
significance. \textit{Journal of the American Statistical Association}, 
\textbf{86}(414), 258--265.

\item  {Hannig, J., Iyer, H., Lai, R. C. S., \& Lee, T. C. M.} (2016).
Generalized fiducial inference: A review and new results. \textit{Journal of
the American Statistical Association}, \textbf{111}(515), 1346--1361.

\item  {Kolmogorov, A. N. } (1933). \textit{Foundations of Probability}.
(Translated by N. Morrison. New York: Chelsea Publishing Company).

\item  {Lee, H. \& Lee, Y.} (2025). Overcoming dilution of collision
probability in satellite conjunction analysis via confidence distribution. 
\textit{Entropy}, \textbf{27}(4), 329.

\item  {Lee, Y.} (2026). Resolving the induction problem: Can we state with
complete confidence via induction that the sun rises forever? To appear at 
\textit{International Statistical Review}.

\item  {Lee, Y. \& Nelder, J. A.} (1996). Hierarchical generalized linear
models. \textit{Journal of the Royal Statistical Society Series B}, \textbf{%
58}(4), 619--656.

\item  {Lee, Y., Nelder, J. A \& Pawitan, Y.} (2017). \textit{Generalized
linear models with random effects: unified analysis via h-likelihood}.
Chapman and Hall.

\item  Neyman, J \& Pearson, E, S, (1933). On the problem of the most
efficient tests of statistical hypothesis. \textit{Philosophical
Transactions of the Royal Society A: Mathematical, Physical and Engineering
Sciences}. \textbf{231}, 289-337.

\item  {Neyman, J.} (1934). On the two different aspects of the
Representative method: the method of stratified sampling and the method of
purposive selection. \textit{Journal of the Royal Statistical Society} B, 
\textbf{97}, 558-625.

\item  {Neyman, J.} (1937). Outline of a theory of statistical estimation
based on the classical theory of probability. \textit{Philosophical
Transactions of the Royal Society A}, \textbf{236}(767), 333--380.

\item  {Pawitan, Y., Lee, H., \& Lee, Y.} (2023). Epistemic confidence in
the observed confidence interval. \textit{Scandinavian Journal of Statistics}%
, \textbf{50}(4), 1859--1883.

\item  {Pawitan, Y. \& Lee, Y.} (2021). Confidence as likelihood. \textit{%
Statistical Science}, \textbf{36}(4), 509--517.

\item  Pawitan, Y. \& Lee, Y. (2024). \textit{Philosophies, puzzles and
paradoxes: a statistician's search for truth}. Chapman and Hall.

\item  {Pitman, E. J. G.} (1957). Statistics and science. \textit{Journal of
the American Statistical Association}, \textbf{52}(279), 322--330.

\item  Reid, N. (2003). Asymptotics and the Theory of inference. \textit{The
Annals of Statistics}. \textbf{31}, 1695-1731.

\item  {Schweder, T. \& Hjort, N. L.} (2016). \textit{Confidence,
likelihood, probability}. Cambridge University Press.

\item  {Skovgaard, I. M.} (2001). Likelihood asymptotics. \textit{%
Scandinavian Journal of Statistics}, \textbf{28}(1), 3--32.

\item  Tang Y. \& Reid N. (2025). Asymptotic behaviour of the modified
likelihood root. \textit{Statistica Sinica}, \textbf{35}, 1559-1581.
\end{itemize}

\color{black}







\end{document}